# The Role of Language in Teaching and Learning Mathematics

*Aslanbek Naziev*
Ryazan State University, Ryazan, Russian Federation
Ulitsa Svobody, 46, Ryazan', Ryazanskaya oblast', Russian Federation, 390000
Tel.: +7 491 297-15-15
a.naziev@rsu.edu.ru

**Abstract**
The purpose of this paper is to emphasize the role of language in the process of teaching and learning mathematics. We will begin with the definition of mathematics given by Cassiodorus (in its essential features repeated in Kolmogorov's famous article "Mathematics" from the "Great Soviet Encyclopedia" and explored in the author's book *Humanitarian oriented mathematics teaching in general school* - Ryazan, Ryazan Institute of Education Development, 2000, 112 p.). This definition illuminates the exceptional role of language in the building, discovering, understanding, teaching, and learning of mathematics. We show that this role is predetermined by the abstract character of mathematics. The abstract object cannot be shown; it can only be described with the help of language. This gives rise to the problem of identifying abstract entities on the grounds of their descriptions. We accentuate the way the language of set theory helps finding the solution of this problem. We illustrate this through the examples from elementary geometry and elementary algebra. And we especially underline the importance of this solution when working with events in the probability theory. After that, we are going to show how proof and language help to form and develop imagery thinking and, thereby, illuminate the possibility of the left brain hemisphere to participate in the formation of imagery thinking.

**1. Introduction**
Most students of mathematics are very disparaging to the language. And they cannot be blamed for this, since the traditional manner of teaching mathematics pushes them to feel this way.
Let us take a look at the lecture on mathematical analysis. The lecturer says:

If functions $f$ and $g$ are continuous at the point $x_0$,
then the function $f+g$ is continuous at the point $x_0$.

But what he writes on the desk is rather different, namely

$$f \quad g \quad x_0 \quad f+g \quad x_0.$$

And this is precisely what his students write in their notebooks, getting accustomed, therefore, to omit words in mathematical texts and implicitly forming in them the belief that mathematics is only the convoluted signs, that words are extraneous, that they are a "literature".
Several years ago, this author invited 24 pupils from a mathematics class of one of the Ryazan schools to solve the following
PROBLEM: Prove that the function defined by the rule "$y = x^2$" is not periodical.
All pupils wrote

$$f(x \pm T) = f(x),$$

and stopped. None of them solved the problem. When the author suggested that they should remember the definition of the periodical function, they said that the boxed equality is nothing but the required definition.

Now, let us remember the genuine
DEFINITION. The function is called periodical if there exists $T \neq 0$ such that for all $x \in D(f)$ we have

$$\boxed{f(x \pm T) = f(x)}.$$

And let us look at the
*Solution of the problem*. Suppose that the function defined by the rule "y=$x^2$" is periodical. Then there exists $T \neq 0$ such that for all

$$x \in \mathbb{R}$$

we have

$$\boxed{(x \pm T)^2 = x^2}.$$

But then, in particular, we have $(T-T)^2 = T^2$, i. e. $0 = T^2$, so $T=0$ — a contradiction. That is all!

As we can see, in order to solve the problem one must only repeat the definition with minor changes. Why were 24 high school mathematics students not able to solve this simple problem? Because they were unable to repeat the definition. But why were they unable to repeat the definition? Because they did not see the words present in the definition. And why is this? Because their teacher told them that in mathematics there should be no words!

To appreciate the role of language in mathematics, one has to start with the definition of mathematics.

**2. What is mathematics?**

Mathematics, unlike other scientific disciplines, can be characterized from two perspectives: from the point of view of the object and on from the point of view of the method.

*2.1. Mathematics from the point of view of matter*

From the point of view of the matter, "mathematics is the science of quantitative relations and spatial forms in the real world.

Quantitative relations (in the general philosophical sense of this term) are characterized, contrary to qualitative, only by their indifference to the specific nature of those objects; they relate...
… spatial forms one may consider as a particular case of quantitative relations.
… mathematics only teaches those relations indifferent to the concrete nature of the objects they relate.
… in the indicated general sense, all relations taught by mathematics always are quantitative." (Kolmogorov (1954); translation into English is mine. — A. N.)

This leads us to the following remarkable definition of mathematics, known as early as VIth century AD (Kolmogorov did not mention it):

> "Mathematica, quam latine possumus dicere 'doctrinalem', scientia est quae abstractam considerat quantitatem. Abstracta enim quantitas dicitur, quam intellectu a materia separantes uel ab aliis accidentibus, ut est par, impar, uel ab alis huiuscemodi, in sola ratiocinatione tractamus."

> "Mathematics is the science of abstract quantity. And abstract quantity is what we learn purely speculatively, distracting it in mind from the substance and random manifestations". (Cassiodorus, (1979), 1· edition VI century AD; (translation into English is mine. — A. N.)

*2.2. Mathematics from the point of view of method*

The abstract character of mathematics predetermines the exceptional role of language in the building, teaching, and learning of mathematics. We cannot show an abstract entity, we can only describe it with the help of language. That is why mathematicians are forced to be extremely exact in their formulations.

## A. MATHEMATICS IS THE PROOF.

"Since the time of the Greeks to say "mathematics" means to say "proof".

(Bourbaki, 1957)

"…mathematics is coextensive with demonstrative reasoning, which pervades the sciences just as far as their concepts are raised to a sufficiently abstract and definite, mathematical-logical level".

(Polya, 1981)

I mention again that the abstract character of mathematics, predetermines the exceptional role of language in the building, teaching, and learning of mathematics. We cannot show an abstract entity, we only can describe it with the help of language. That is why mathematicians are forced to be extremely exact in their formulations.

We cannot show a geometrical point because the point has no dimensions. We cannot show a geometrical line, firstly because the line has no thickness, and secondly, because it is infinite in both directions. That is why we are forced to describe points, lines, and the relations among them with the help of axioms, and, in the last analysis, with the help of language. And we must understand that "points", "lines", "triangles" etc. from pictures are not points, lines, triangles themselves, but only their imperfect representations.

We cannot show numbers. For example, in the following picture, we can see three cats. Well, we see the *cats*, but where is *three*?

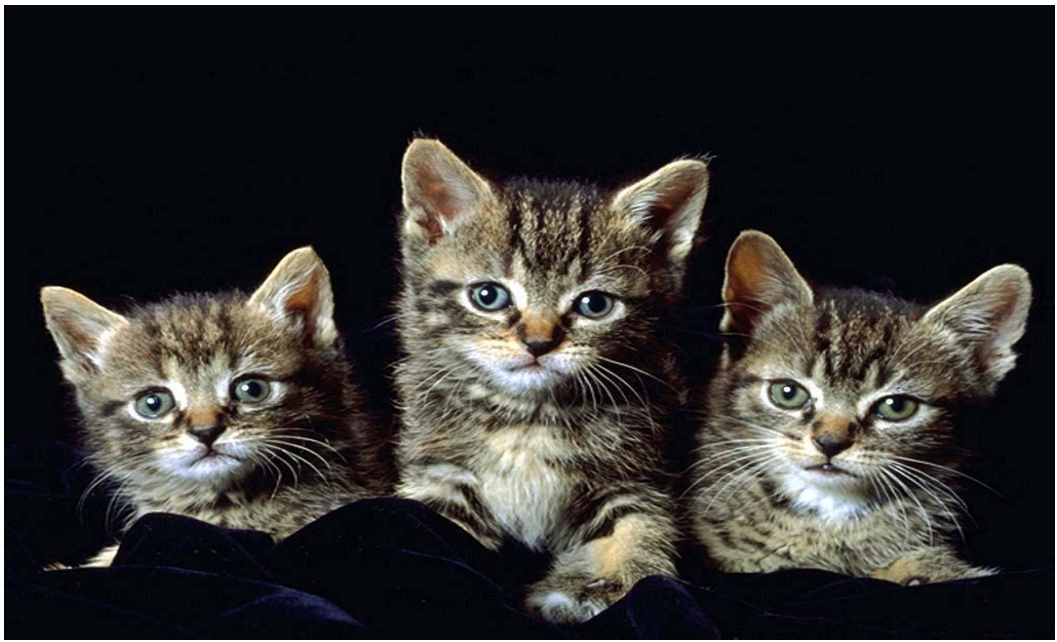

*Figure 1. Three cats*

*Three* is connected with this picture mentally. *Three* cannot be shown, it only may be thinkable.

Oh, maybe this

3,

is number three? No! Why? Because, if we agree that this is number three, then we must agree that this

# 3

is also number three but the second number is greater than the first one! So, we must distinguish between numbers, abstract entities, and their representations with the help (of signs) of language.

Analogously, in probability theory, we cannot show events. Let us imagine that we roll a dice and it falls down as is shown in the following picture (Figure 2a):

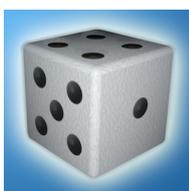 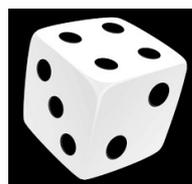

*Figure 2a. A dice*     *Figure 2b. A dice*

What event can we see? None! Why? Because events, similarly to all other mathematical objects, are *abstract* entities. (If you think that you see *the events*, try to explain what events do you see on Figs 2a and 2b, identical or different?) We do not (and cannot!) see *the event*, we see *a dice* with four points on the top face, five points on the right face and one point on the front face. Events are connected with this picture mentally, and we have a large amount of such events:
- it is rolled a dice with four points,
- it is rolled a dice with an even number of points,
- it is not rolled a dice with less than 4 points (3,2,1,0),
- it is not rolled a dice with more than 4 points (5,6,7 etc).

There is an infinite number of events that can be mentally associated with the above single picture and none of them can be mentioned as *the event* represented by this picture.

Now, the question arises naturally: how can we learn about these events? The answer is obvious from the preceding paragraph: we learn that with the help of language, namely, with the help of linguistic descriptions of these abstract entities.

After we realize this, we will face the next problem: how to identify or distinguish the events on the grounds of their descriptions?

Similar questions arise in all other branches of mathematics. Are the properties of real numbers expressed by the sentences "$x < 0$" and "$x^3 < 0$" distinct or identical? Are the geometric figures described as "the locus of points equidistant from the ends of the line segment" and as "the line through the middle of the line segment perpendicular to the segment" distinct or identical? Are the functions defined by the conditions "$y = \cos 2x$" and "$y = \cos^4 x - \sin^4 x$" distinct or identical? Etc.

The solution of all of these problems was found by Georg Cantor (1845–1918). This solution is contained in the rule: translate your problems to the language of set theory where you have exact criteria of identity and difference. Or, as Quine (1966) wrote,

> "One of the sources of clarity in mathematics is the tendency to talk of classes rather than properties. Whatever is accomplished by referring to a property can generally be accomplished at least as well by referring to the class of all the things that have the property. Clarity is gained in that for classes we have a clear idea of sameness and difference: it is a question simply of their having the same or different members."

For example, instead of thinking about *properties*, expressed by the sentences "$x < 0$" and "$x^3 < 0$" we should instead think about the *sets* of real numbers $x$ such that $x < 0$, $x^3 < 0$ respectively. By doing so, we can see that these two *sets* are identical and we can express this by saying that the *properties* are identical. Therefore, due to Cantor, mysterious properties became fully intelligible sets and the problems regarding their difference or identity vanished.

Similarly, the locus of points equidistant from the ends of the line segment and the line through the middle of the line segment perpendicular to the segment became identical sets of points and the problem vanishes.

Again, the functions defined by the conditions "$y = \cos 2x$" and "$y = \cos^4 x - \sin^4 x$" became identical sets of ordered pairs of real numbers, and the problem also vanishes.

Now turn to the probability theory. Here is the solution as it follows: first of all, divide all events into elementary and compound events. Elementary events are those which cannot be decomposed into the other events, compound events are those which can be. This may be expressed in other words. Elementary events are those which come directly, compound events are those which come indirectly, by means of other events. Elementary events are called outcomes or samples.

For example, when we roll a dice, the event consisting of dropping out of six points comes directly, whereas the event consisting in the dropping of an even number of points does not come directly, it only comes by means of other events: the dropping of two, of four, of six points. Thus, all events, elementary or not, come by means of elementary events, or samples. Now, the solution in the spirit of Cantor of the problem of difference and identity of events consists in set-theoretic interpretation of events. According to this interpretation, every event is seen as the set of samples in which given event is dropped. Due to this principle, we get exact rules for deciding all questions regarding the differentness and identicalness of events.

For example, let us toss a dice and consider these events:
- more than zero points dropped;
- less than seven points dropped.

*Descriptions* of these events are quite different. And what could be said about *the events*? The *sets of samples* corresponding to these descriptions are identical, and this enables us to conclude (according to the set-theoretic interpretation of events) that these *events* are identical.

### B. MATHEMATICS TEACHING AND FORMATION OF IMAGERY THINKING

Are we able to see with our eyes a line and a circle which have exactly one common point? No! We can draw a line and a circle a thousand times, and every time we will see either no common point, two common points, or the whole segment of common points (look at Figure 3 below).

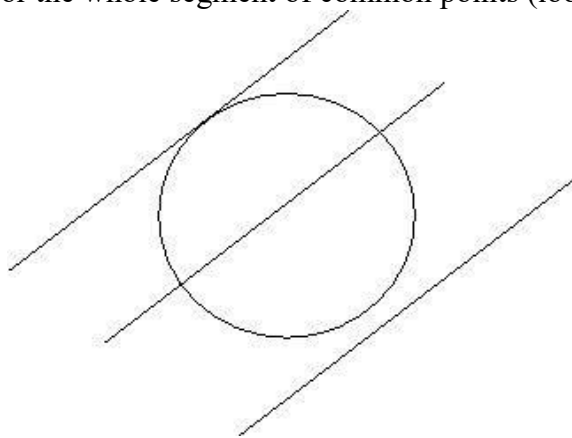

*Figure 3. Circle and Lines*

So it is impossible to see that with our own eyes. We only can "see" this in our mind with the help of proof. ("Not with the eyes of the head, but with the eyes of the mind", as Galilei liked to say.) Therefore, if we restrict ourselves to looking at pictures, we will have an incorrect image. To have the correct image, we need a proof. And we can give it.

*Proof.* Let $\omega$ be a circle with the center $O$ and radius $R$. Take a point $A$ on the circle. Draw segment $OA$ and the line $l$ passing through $A$ perpendicular to $OA$. Then $OA$ is the perpendicular dropped to $l$ from $O$ (and its length equals $R$). And for any other point $B$ of the line $l$, the segment $OB$ is inclined to the line $l$ and so $OB > OA = R$, i. e. $\omega$ and $l$ has only one common point (namely, $A$).

This example, among other things, shows how the image is forming in a verbal-logical way.

Looking at pictures creates the impression that the line and the circle cannot have only one common point. But this is only the *beginning* of the process of reflection of the studied phenomena in our mind, while the image is the *result* of reflection. And if we stop at the stage of looking, we will have a wrong image. Reasoning shows us the fallacy of representations based on direct perception and leads to an adequate reflection of the phenomenon in our mind. As a result, the right image is formed with the help of reasoning, that is, in a verbally-logical way.

Thus, the proof is a very important tool for the formation of imagery thinking (contrary to looking!).

Let us illustrate this with the help of another

EXAMPLE. Four spheres of the same radius touch each other pair-wise externally. Find the radius of the sphere that touches each of these spheres internally.

The author repeatedly gave this task to the graduates of the school. More than 40% of the decision-makers didn't solve this problem. They have arbitrarily replaced the given spheres with the circles and began from the picture represented in Figure 4a, from which it is clear that they misunderstood the condition "touch each other pair-wise". Instead of the right "each touches each" they interpreted this condition in their own way, believing that the circles are divided into pairs touching each other. (A similar error occurs when one says that a parallelogram is a quadrilateral whose sides are pair-wise parallel.) They also replaced the fifth sphere with a circle and at the same time incorrectly interpreted the condition "touches internally". Instead of remembering what this condition means (each pair of touched spheres have a common tangent plane and lie on one side of it), they believed that it is talked about a sphere enclosed in the inner area bounded by these spheres and found the radius of a small circle represented on Figure 4b (not noticing that the case they are considering refers to an external touch).

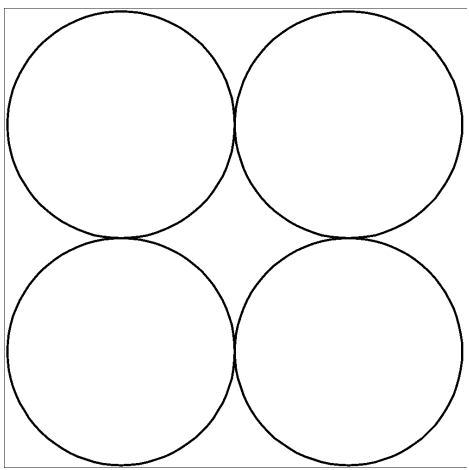   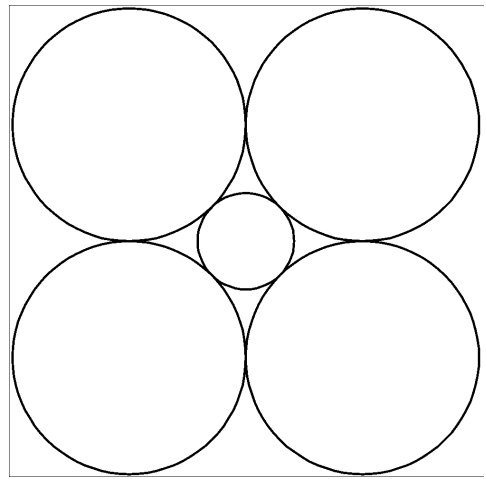

*Figure 4a. Four circles*          *Figure 4b. Five circles*

Two gross errors made due to a misunderstanding of the conditions of the problem (and therefore of a linguistic nature). It is clear that in this case there can be no question of forming the right images.

Now we are going to create the right images. Let's start with four equal circles since this error was made. We want to form the correct image of four circles of the same radius on the plane, pair-wise touching each other externally. Is it possible, and if possible, how?

ANSWER: this is impossible. *Indeed*, if two circles on a plane touched externally, the distance between their centers is equal to the sum of their radii. Let the radii of all four circles be $r$. Consider three of these circles. In view of what has just been said, their centers are located at the vertices of a regular triangle with the side of length $2r$. The center of the fourth circle is equidistant from the centers of these three circles and, hence, is the center of the circle circumscribed near the regular triangle with a side of length $2r$. But, on the other hand, it must be away from each of these centers by $2r$. It turns out that the radius of a circle, which is circumscribed near a regular triangle with the side of length $2r$, is equal to $2r$, which is incorrect. Hence, on the plane, four circles of the same radius cannot pair-wise touch each other.

Another example one can find is in the domain of spheres. There, it is possible to find four spheres of equal radii pairwise touching each other and in one way only, namely, in the way where the centers of the spheres are located at the vertices of a regular tetrahedron with the side of length $2r$. In order to obtain a sphere touching the four given spheres internally, we take the sphere circumscribed near this tetrahedron, enlarge its radius by $r$, and prove that obtained sphere is the required one.

So, we have found the required sphere. Let us stress again that that was done by means of the language and proof, that is, in a verbally-logical way.

Let us consider a more complicated

EXAMPLE. How many roots has the equation
$$\left(\frac{1}{16}\right)^x = \log_{\frac{1}{16}} x?$$

In order to answer this question, note that the functions
$$y = \left(\frac{1}{16}\right)^x \text{ and } y = \log_{\frac{1}{16}} x$$
are mutually inverse, so their graphs are symmetrical by the bisector of the first and third coordinate angles, and because each of them intersects this bisector, they intersect it in one point. Besides, these functions are both decreasing, and it seems that their graphs go something like this:

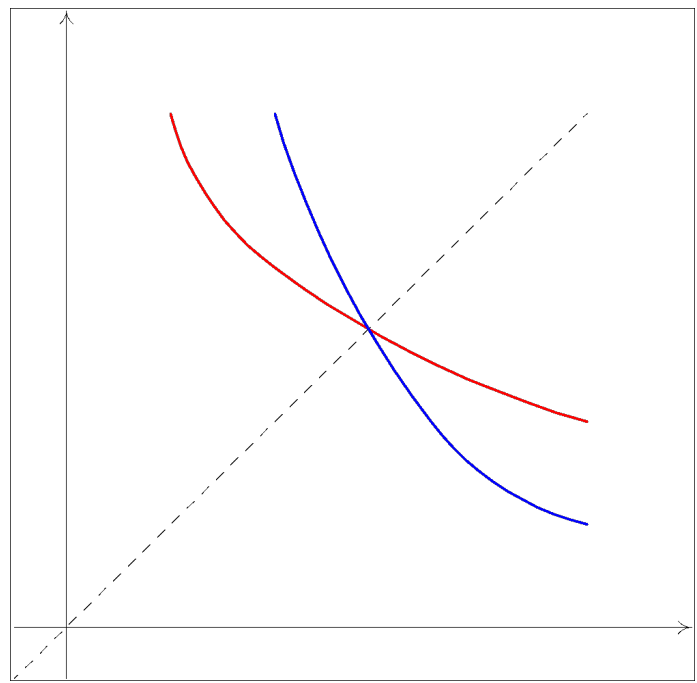

*Figure 5. "Evident" behavior of the graphs of the functions* $y = \left(\frac{1}{16}\right)^x$ *and* $y = \log_{\frac{1}{16}} x$.

It seems obvious that the graphs of our functions have only one common point, and thus the given equation — only one root.

But wait! We see other roots, $\frac{1}{2}$ and $\frac{1}{4}$! So, how many roots does this equation really have? We tried asking a computer to give the answer. The computer gave the following answer: Figure 6, where in blue is the graph of the function $y = \left(\frac{1}{16}\right)^x$, and in red is the graph of the function $y = \log_{\frac{1}{16}} x$.

This entire story shows us again that we cannot rely on the evidence. We need proof!

Thorough analysis with the help of the first and second derivatives (Vilenkin, 1980) proves that given equation has exactly three roots, $\frac{1}{2}$, $\frac{1}{4}$, and one more root corresponding to the point lying on the mentioned bisector; this root cannot be found by elementary methods. Graphs of the functions go like the blue and the red curves in Figure 7, only much closer one to another.

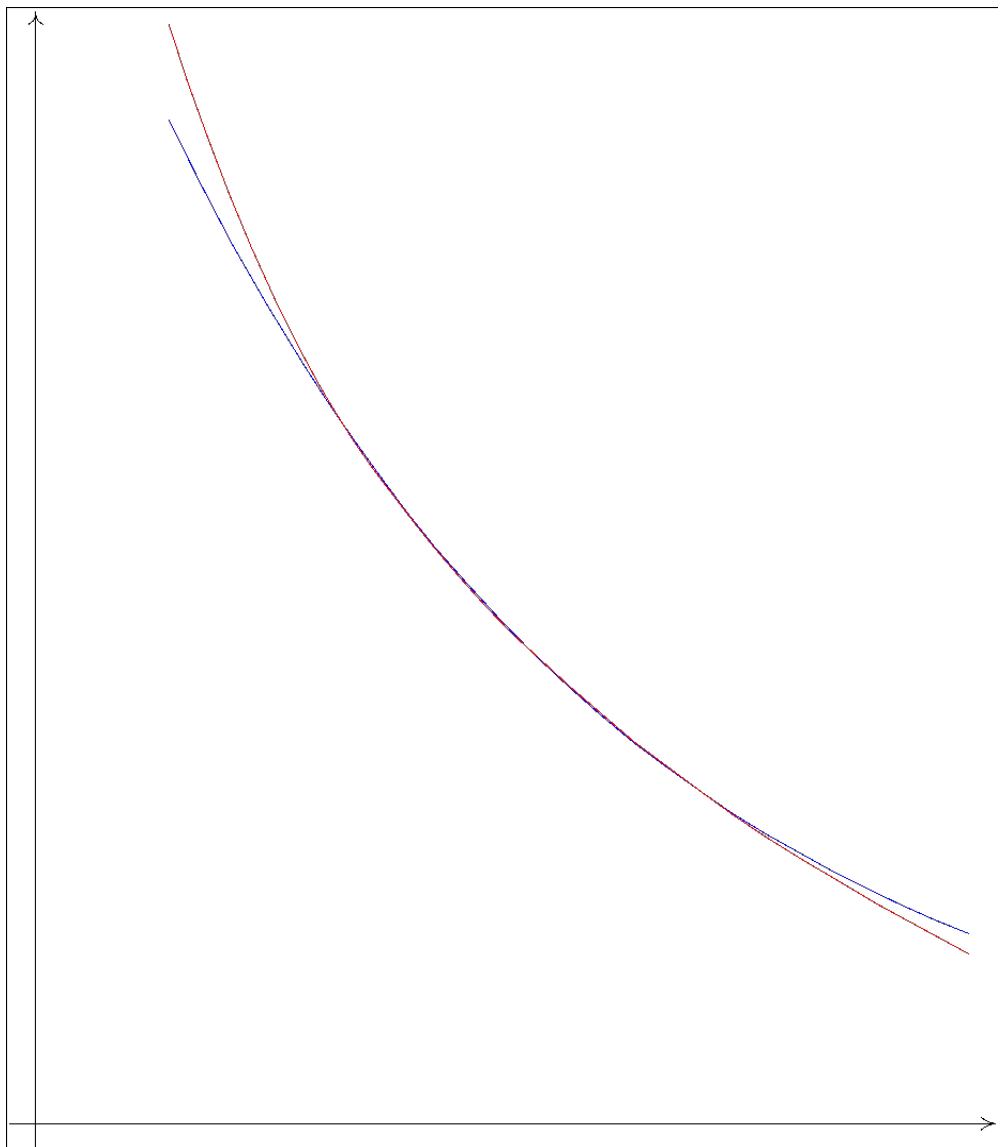

*Figure 6. Computer's answer*

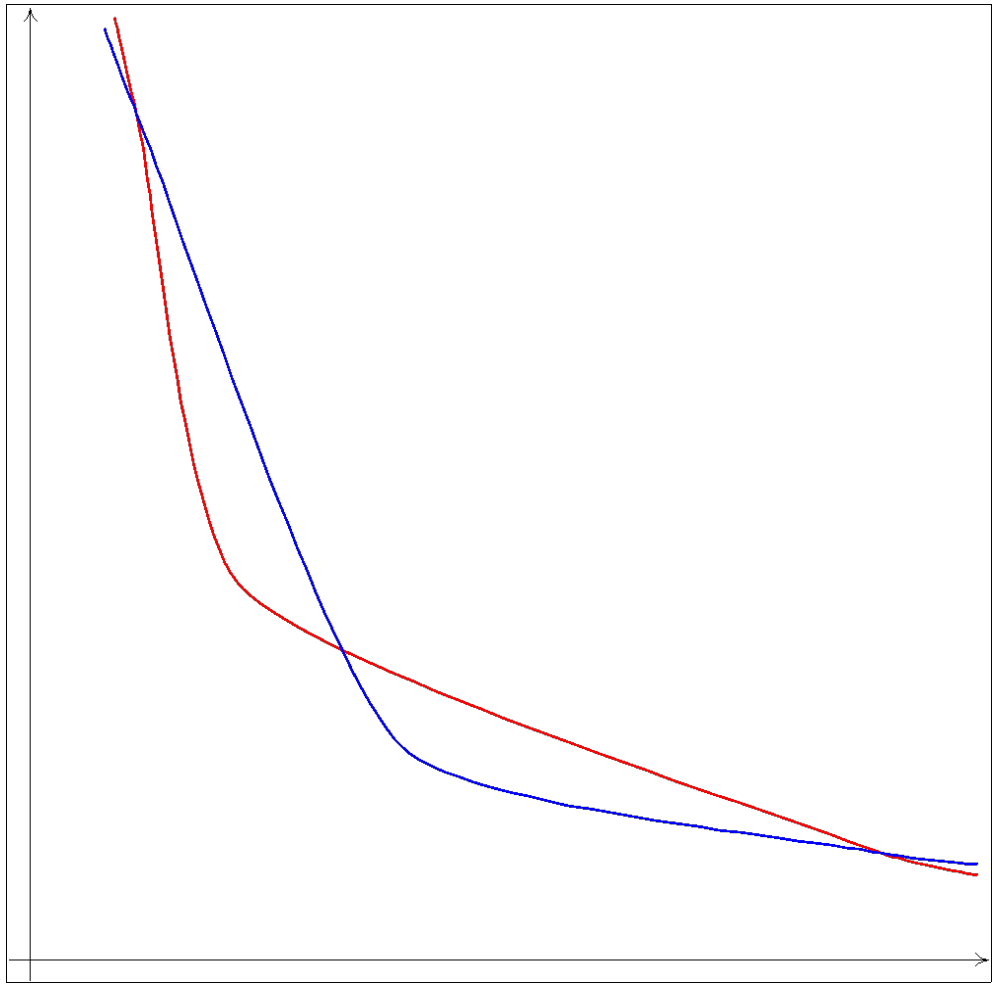

*Figure 7. Schematic drawing of the graphs of the functions after research*

This is the right image of this situation; not precise, but right, in the sense that it rightly reflects the joint behavior of the graphs of the functions. It is important to mention, that, again, the right image is obtained due to proof.

All the above examples show that proof is an indispensable tool for obtaining the right images. Now, note that it couldn't be otherwise! Remember what A. N. Whitehead wrote about proof, understanding, and self-evidence:

> … understanding is not primarily based on inference. Understanding is self-evidence. … But the large self-evidence of mathematical science is denied to humans. …
> …our clarity of intuition is limited, and it flickers. Thus inference enters as means for the attainment of such understanding as we can achieve. Proofs are the tools for the extension of our imperfect self-evidence.

### 3. Language, proof, and functional brain asymmetry

The researcher of the processes of intelligence and creativity Paul Torrens (1915–2003) was among the first who paid attention to the working features of the hemispheres from a person's brain. The scientist conducted an experiment in the course of which four types of thinking were established:
- left hemisphere thinking, as one that is built on logic and analysis;
- right hemisphere thinking, where the thinking process is directed by emotions, intuition, and images;
- mixed thinking, where both right and left hemispheres are equally active, each of which is switched on at the right moment;

- integrated thinking, when the right hemisphere and left hemisphere thinking work simultaneously.

Torrens stressed that among these types of thinking there are no good or bad types: everyone has advantages and disadvantages.
However, the scientists and educators in recent years have increasingly attached importance to the development of right hemisphere thinking.

They rightly point out that the classical education system is built in such a way as to develop predominantly left hemisphere thinking, almost completely ignoring the development of right hemispheric skills, and the main culprit of this state of affairs they believe is mathematics. But this writer is sure that it's not the mathematics to blame, but the traditional system of teaching, including mathematics teaching.

The author suggested the conception of humanitarian oriented mathematics teaching (Naziev, 2000) according to which "to teach mathematics means to systematically encourage the students to discover their own proofs". He was repeatedly criticized for this position, seeing in it an emphasis on verbal and logical at the expense of figurative. All the examples considered above are made in the spirit of this conception. And all these examples demonstrate that with a correct understanding of the verbally-logical, not only that it does not damage the figurative, but on the contrary, it contributes to its development. This shows that the danger of left-hemispheric thinking is exaggerated. But the danger of the right hemispheric is very real. This danger (without mentioning it explicitly) warned (Huizinga, 1936, in inverse translation from Russian):
"It is the increased visual suggestibility (read — right-hemisphereness) that Achilles heel, which beats modern people advertisement, using the weakening ability of judgment, the ability to think and evaluate independently. This applies to both commercial and political advertisement."

A person permeated by a holistic perception loses the ability to analyze and critically evaluate the information coming to him. And as a result, he is subject to the thoughtless repetition of thoughts and phrases imposed on him from outside. As a result, the desire to defeat verbalism generates the worst version of it, for which Huizinga finds an extremely accurate name, verbicide.

One artist, familiar to the author, likes to repeat: do not go deep, do not analyze, look at it as a whole, take it in one piece. This type of artistic perception and behavior is especially appealing to despotic bosses who want their orders to be executed without hesitation, without delving into them or analyzing. So, the seemingly good aspiration to the integrity of perception opens the way to conformism and through it to tyranny.

Those who are inclined to overvalue the creative potential of the right hemisphere should carefully think over the following words by (Ivanov, 1994); just one phrase full of profound meaning: "The right hemisphere in the norm (when censorship of the left hemisphere is carried out) is wordless, and this is the source of his creative potential." Therefore, to uncover the creative potential of the right hemisphere, censorship of the left hemisphere is necessary. And what happens when the right hemisphere is freed from this censorship is written one page earlier: "self-killing is the limiting case that can be characterized as the killing of the left hemisphere by the right hemisphere."

### 3. Conclusions

Our research has shown that teaching mathematics should be done in a humanitarian way. Teaching mathematics has to be an activity directed to the discovering of proofs, that connects the left hemisphere to help the right hemisphere to build the right images, thereby, opening the way to "achieving a perfect harmony of the left and right", the necessity of which was repeatedly stressed by the outstanding Russian neurolinguist Vyach Vs. Ivanov. Moreover, although it seems strange, teaching mathematics in the spirit of the author's conception, contributes, albeit small, to the prevention of tyranny and suicide.

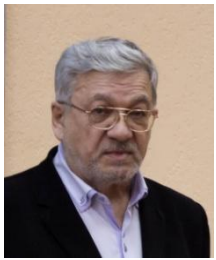

**Aslanbek NAZIEV** (b. April 3 1946) received his master's degree in Mathematics and Programming, Ph D degree (Russian "Candidate of Sciences") in Mathematics and Sc D degree (Russian "Doctor of Sciences") in Mathematics Education, all from the Moscow State Pedagogical University. Now he is full professor in department of Mathematics and Mathematics Teaching, Faculty of Physics and Mathematics, Ryazan State University, Ryazan, Russian Federation. His current research interests include Philosophy of (Mathematics) Education, Humanitarian Oriented Mathematics Teaching, Logical problems of Mathematics Education, Integrating of ICT and Mathematics in Mathematics and Science Education, Multidimensional Education, and many other topics. He has (co-)authored more than 120 printed works (including 7 books) and more than 100 electronic publications, participated a large amount of scientific conferences (13 just in 2017), in many conferences he was a key-note speaker (4 just in 2017).